\documentclass{amsart}
\usepackage{amsmath}
\usepackage{amssymb} 
\usepackage{amscd} 
\usepackage{graphicx} 
\usepackage{epsfig}
\usepackage[dvips]{color}

\newtheorem{theorem}{Theorem}[section]
\newtheorem{lemma}[theorem]{Lemma}

\newtheorem{proposition}[theorem]{Proposition}

\newtheorem{claim}[theorem]{Claim}

\newtheorem{conjecture}[theorem]{Conjecture}
\newtheorem{remark}[theorem]{Remark}

\theoremstyle{definition}
\newtheorem{definition}[theorem]{Definition}

\numberwithin{equation}{section}
\numberwithin{figure}{section}
\numberwithin{table}{section}

\begin{document}
\baselineskip 14pt

\title{The slope conjecture for graph knots}

\author[K. Motegi]{Kimihiko Motegi}
\address{Department of Mathematics, Nihon University, 
3-25-40 Sakurajosui, Setagaya-ku, 
Tokyo 156--8550, Japan}
\email{motegi@math.chs.nihon-u.ac.jp}

\author[T. Takata]{Toshie Takata}
\address{Graduate School of Mathematics, Kyushu University, 
744 Motooka, Nishi-ku, Fukuoka 819--0395, Japan}
\email{ttakata@math.kyushu-u.ac.jp}

\dedicatory{}

\begin{abstract}
The slope conjecture proposed by Garoufalidis asserts that the Jones slopes given by the sequence of degrees of the colored Jones polynomials are boundary slopes. 
We verify the slope conjecture for graph knots, i.e. knots whose Gromov volume vanish. 
\end{abstract}

\maketitle

{
\renewcommand{\thefootnote}{}
\footnotetext{
}

\section{Introduction}
\label{section:Introduction}

Let $K$ be a knot in the $3$--sphere $S^3$ and $E(K)$ the exterior $S^3 - \mathrm{int}N(K)$. 
Denote by $(\mu, \lambda)$ the preferred meridian-longitude pair of $K$. 
Then any homotopically nontrivial simple closed curves in $\partial E(K)$ 
represents $p[\mu] + q[\lambda] \in H_1(\partial E(K))$ for some relatively prime integers $p$ and $q$. 
We call $p/q \in \mathbb{Q} \cup \{ \infty \}$ a \textit{boundary slope} of $K$ if there exists a connected, orientable, incompressible and 
boundary-incompressible surface $F$ such that a component of $\partial F$ represents $p[\mu] + q[\lambda] \in H_1(\partial E(K))$. 
Let us define: 
$$bs(K) = \Bigl\{ r \in \mathbb{Q}\cup \{\infty\}\ |\ r\ \textrm{is a boundary slope of}\ K\Bigr\}.$$
Following Hatcher \cite{Hat1} $bs(K)$ is a finite subset of $\mathbb{Q} \cup \{\infty\}$ for every knot $K$. 

\smallskip

The \textit{colored Jones function} of $K$ is a sequence of Laurent polynomials $J_{K, n}(q) \in \mathbb{Z}[q^{\pm 1}]$ for $n \in \mathbb{N}$, 
where $J_{K, 2}(q)$ is the ordinary Jones polynomial of $K$. 
Let $\delta_K(n)$  be the maximum degree of $J_{K, n}(q) \in \mathbb{Z}[q^{\pm 1}]$. 
We call $x \in \mathbb{R}$ a \textit{cluster point} of a sequence $\{ x_n \}$ if $x$ is a limit point of a subsequence of $\{ x_n \}$. 
We define $js(K)$ as follows: 
$$js(K) = \Bigl\{ \textrm{cluster points of the sequence}\ \{\frac{4 \delta_K(n)}{n^2} \}_{n \in \mathbb{N}} \Bigr\}.$$
Since the colored Jones function is $q$--holonomic \cite[Theorem~1]{GL}, 
Theorem~1 in \cite{Garoufalidis} shows $\delta_K(n)$ is a \textit{quadratic quasi-polynomial}, i.e. 
$$
\delta_K(n)=c_2 (n) n^2 + c_1 (n) n+c_0(n)
$$
for rational valued periodic functions $c_i(n)$ with an integral period. 
By Lemma~1.8 in \cite{Garoufalidis}, $js(K)$ is the finite set of 
4 times the rational values  of the periodic function $c_2 (n)$. 
Using the minimum degree $\delta^*_K(q)$ of $J_{K, n}(q)$ instead of $\delta_K(n)$, 
we can define: 
$$js^*(K) = \Bigl\{ \textrm{cluster points of the sequence}\ \{\frac{4 \delta^*_K(n)}{n^2} \}_{n \in \mathbb{N}} \Bigr\}.$$
As noted in \cite[1.4]{Garoufalidis}, 
$\delta^*_K(n) = -\delta_{K^*}(n)$ and thus 
$js^*(K) = -js(K^*)$, 
where $K^*$ is the mirror image of $K$ and $-X := \{ -x_1,\dots, -x_m\}$ if $X = \{ x_1,\dots, x_m\}$. 
We call an element in $js(K) \cup js^*(K)$ a \textit{Jones slope} of $K$. 

In \cite{Garoufalidis}, Garoufalidis proposed the following conjecture which relates Jones slopes and boundary slopes. 

\begin{conjecture}[\textbf{Slope conjecture}]
\label{slope conjecture}
For any knot $K$, 
every Jones slope is a boundary slope, 
i.e. $js(K) \cup js^*(K) \subset bs(K)$. 
\end{conjecture}

The conjecture was verified for torus knots, 
some non-alternating knots,  the $(-2, 3, p)$--pretzel knots \cite{Garoufalidis}, 
adequate knots \cite{FKP} and a $2$--parameter family of $2$--fusion knots \cite{DG, GV}. 
Note that the class of adequate knots includes all alternating knots and most Montesinos knots. 
Recently, the conjecture was verified for iterated cables of adequate knots and for iterated torus knots. \cite{KT, KT2}.

In the present note we give further supporting evidence for the sloe conjecture by proving: 

\begin{theorem}
\label{slope conjecture composite}
Suppose that $K_1$ and $K_2$ satisfy the slope conjecture, 
then the connected sum $K_1 \sharp K_2$ also satisfies the slope conjecture. 
\end{theorem}

A knot $K$ is called a \textit{graph knot} if its exterior $E(K)$ is a graph manifold, 
i.e. there is a family of tori which decomposes $E(K)$ into Seifert fiber spaces. 
This implies that any graph knot is obtained from unknots by a finite sequence of operations of cabling and connected sum; 
see \cite[Corollary~4.2]{Gor}. 
A graph knot can be also characterized as a knot whose Gromov volume vanishes \cite{Gr, So, Th}. 

As an application of Theorem~\ref{slope conjecture composite} and \cite{KT} we establish: 

\begin{theorem}
\label{slope conjecture graph}
Every graph knot satisfies the slope conjecture. 
\end{theorem}

\bigskip

\noindent
\textbf{Acknowledgments.}
We would like to thank Effie Kalfagianni and Anh Tran for suggesting an error in an earlier version of the paper. 

K.M. has been partially supported by JSPS Grants--in--Aid for Scientific 
Research (C), 26400099, The Ministry of Education, Culture, Sports, Science and Technology, Japan and Joint Research Grant of Institute of Natural Sciences at Nihon University for 2014. 
T.T. has been partially supported by JSPS Grants--in--Aid for Scientific 
Research (C), 25400094, The Ministry of Education, Culture, Sports, Science and Technology, Japan.

\section{The slope conjecture and connected sum operation}
\label{composite knot}

In this section we prove Theorem~\ref{slope conjecture composite}, 
i.e. the positivity of the slope conjecture is preserved under connected sum. 

\medskip

\noindent
\textit{Proof of Theorem~\ref{slope conjecture composite}.}
First we describe jones slopes $js(K_1 \sharp K_2)$.

\begin{lemma}
\label{jones slope composite}
If $p/q \in js(K_1 \sharp K_2)$, 
then there exist $p_1/q_1 \in js(K_1)$ and $p_2/q_2 \in js(K_2)$ such that $p/q = p_1/q_1 + p_2/q_2$. 
\end{lemma}

\noindent
\textit{Proof of Lemma~\ref{jones slope composite}.}
Let us put 
$\delta_{K_i}(n)=\alpha_i(n)n^2+\beta_i(n)n+\gamma_i(n)$ for $i=1,2$, 
where $\alpha_i(n)$, $\beta_i(n)$ and $\gamma_i(n)$ are periodic functions with integral periods. 
It is known that 
$$
 J_{K_1 \sharp K_2,n}(q)=J_{K_1,n}(q)J_{K_2,n}(q)
$$ 
and so we have: 
\begin{eqnarray*}
\delta_{K_1 \sharp K_2}(n)&=&\delta_{K_1}(n)+\delta_{K_2}(n) \\
 & & (\alpha_1(n)+\alpha_2(n))n^2+(\beta_1(n)+\beta_2(n))n+(\gamma_1(n)+\gamma_2(n)).
\end{eqnarray*}
Since $\alpha_i(n)$, $\beta_i(n)$ and $\gamma_i(n)$ are periodic functions with integral periods, 
so are $\alpha_1(n)+\alpha_2(n)$, $\beta_1(n)+\beta_2(n)$, 
and $\gamma_1(n)+\gamma_2(n)$. 
Therefore, the Jones slope $p/q$ of $K_1 \sharp K_2$ is an element of the finite set of 
the rational values of $4(\alpha_1(n)+\alpha_2(n)) = 4\alpha_1(n) + 4\alpha_2(n)$, 
and hence $p/q = p_1/q_1 +p_2/q_2$ for some Jones slopes $p_1/q_1 \in js(K_1), p_2/q_2 \in js(K_2)$. 
\hspace*{\fill} $\square$(Lemma~\ref{jones slope composite})

\medskip

Since every composite knot has an essential meridional annulus, 
in the following we may assume $q_i > 0$ for $i = 1, 2$. 

\begin{lemma}
\label{boundary slope composite}
If $p_1/q_1 \in bs(K_1)$ and $p_2/q_2 \in bs(K_2)$, 
then $p_1/q_1 + p_2/q_2 \in bs(K_1 \sharp K_2)$. 
\end{lemma}

\noindent
\textit{Proof of Lemma~\ref{boundary slope composite}.}
Let $A$ be an essential annulus in $E(K_1 \sharp K_2)$ which decomposes $E(K_1 \sharp K_2)$ into 
$E(K_1)$ and $E(K_2)$. 
Let $F_i$ be an essential surface in $E(K_i)$ and $m_i$ the number of boundary components of $F_i$. 
(If $p_i \ne 0$, then for homological reason $m_i$ is an even integer.)
Then $\partial F_i$ consists of $m_i$ mutually parallel loops each of which has slope $p_i/q_i$ $(q_i > 0)$. 
Note that the core of $A$ is a meridian of $K_i$ and choose $F_i$ so that $A \cap F_i$ consists of $m_iq_i$ spanning arcs in $A$. 
See Figure~\ref{composite1}. 

\begin{figure}[!ht]
\begin{center}
\includegraphics[width=0.6\linewidth]{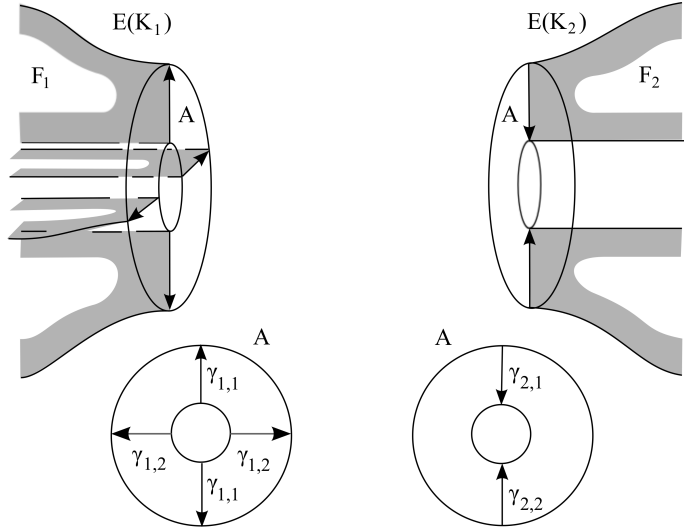}
\caption{Essential surfaces $F_1 \subset E(K_1)$ and $F_2 \subset E(K_2)$ with $m_1 = 2,\ q_1 = 2,\ m_2 = 2,\ q_2 = 1$; 
$\partial F_1 = \gamma_{1,1} \cup \gamma_{1,2},\ \partial F_2 = \gamma_{2, 1} \cup \gamma_{2,2}.$}
\label{composite1}
\end{center}
\end{figure}

Orient $\partial F_i$ so that they run the same direction on $\partial E(K_i)$ (independent of an orientation induced from $F_i$) and 
a component of $\partial F_1$ and that of $\partial F_2$ has opposite orientations on $A$ as in Figure~\ref{composite1}. 
In general, 
$m_1q_1 \ne m_2q_2$, 
i.e. the number of components of $A \cap F_1$ does not coincide with that of $A \cap F_2$, 
so we take $m_2q_2$ parallel copies of $F_1$ and $m_1q_1$ parallel copies of $F_2$. 
Let us denote these (disconnected) surfaces by $m_2q_2F_1 \subset E(K_1)$ and $m_1q_1F_2 \subset E(K_2)$, respectively. 
We give an orientation on the boundary of $m_2q_2F_1$ (resp. $m_1q_1F_2$) 
so that it coincides with that of $\partial F_1$ (resp. $\partial F_2$). 
Since both $A \cap m_2q_2F_1$ and $A \cap m_1q_1F_2$ consist of 
$m_1m_2q_1q_2$ spanning arcs in $A$, 
we can connect $m_2q_2F_1$ and $m_1q_1F_2$ along the annulus $A$ to obtain a possibly disconnected surface $F'$ 
in $E(K_1 \sharp K_2)$. 
Note that all the components of $\partial F'$ run the same direction on  $\partial E(K_1 \sharp K_2)$ 
with respect to the orientation given in the above.

\begin{claim}
\label{slope F'}
Each component of $F' \cap \partial E(K_1 \sharp K_2)$ has slope $p_1/q_1 + p_2/q_2$. 
\end{claim}

\noindent
\textit{Proof of Claim~\ref{slope F'}.}
Let $(\mu_i, \lambda_i)$ and $(\mu, \lambda)$ 
be preferred meridian-longitude pairs of $K_i$ and $K_1 \sharp K_2$; 
we take $\mu_1 = \mu_2  = \mu \subset \partial A$. 
Orient them so that $\langle \mu_i, \lambda_i \rangle = \langle \mu, \lambda \rangle = 1$ and 
$\langle \mu, \partial F' \rangle  
= \langle \mu_1 , \partial (m_2q_2F_1) \rangle 
= \langle \mu_2 , \partial (m_1q_1F_2) \rangle > 0$, 
where $\langle \alpha, \beta \rangle$ denotes the algebraic intersection number between $\alpha$ and $\beta$. 
Then 
$$\langle \mu, \partial F' \rangle 
= \langle \mu_1, \partial(m_2q_2F_1) \rangle 
= m_2q_2 \langle \mu_1, \partial F_1 \rangle = m_2q_2(m_1q_1) = m_1m_2q_1q_2 $$ 
and 
\begin{eqnarray*}
\langle \partial F', \lambda \rangle &=& 
\langle \partial(m_2q_2F_1), \lambda_1 \rangle + \langle \partial(m_1q_1F_2), \lambda_2 \rangle \\
&=& m_2q_2\langle \partial F_1, \lambda_1 \rangle + m_1q_1 \langle \partial F_2, \lambda_2 \rangle \\
&=& m_2q_2(m_1p_1) + m_1q_1(m_2p_2) \\
& =&m_1m_2p_1q_2 + m_1m_2q_1p_2.
\end{eqnarray*}
Thus $F' \cap \partial E(K_1 \sharp K_2)$ represents  
\begin{eqnarray*}
 & & (m_1m_2p_1q_2 + m_1m_2q_1p_2)[\mu] + m_1m_2q_1q_2 [\lambda] \\
 & & = m_1m_2((p_1q_2+q_1p_2)[\mu] + q_1q_2[\lambda]) \in H_1(\partial E(K_1 \sharp K_2)).
\end{eqnarray*} 
Let $k$ be the greatest common divisor of $p_1q_2+q_1p_2$ and $q_1q_2$. 
Then $F' \cap \partial E(K_1 \sharp K_2)$ consists of $m_1m_2k$ parallel loops each of which 
has slope $$\frac{(p_1q_2+q_1p_2)/k}{(q_1q_2)/k} = (p_1q_2+q_1p_2)/q_1q_2 = p_1/q_1 + p_2/q_2.$$  
\hspace*{\fill} $\square$(Claim~\ref{slope F'})

\medskip

Let $F$ be a connected component of $F'$. 
If $F$ is non-orientable, 
then we take a tubular neighborhood $N(F)$ of $F$ in $E(K_1 \sharp K_2)$ and 
we replace $F$ by $\partial N(F)$, 
which is an orientable double cover of $F$ and each component of $\partial N(F)$ has slope $p_1/q_1 + p_2/q_2$; 
for simplicity we continue to use the same symbol $F$ to denote $\partial N(F)$. 
Since $F_i \subset E(K_i)$ is orientable, 
$F \cap E(K_i)$ consists of parallel copies of $F_i$ for $i = 1, 2$. 
Note also that for each component of $F \cap E(K_i)$, 
its boundary component across $A$ in the same direction. 

\begin{claim}
\label{essential}
The surface $F$ is essential in $E(K_1 \sharp K_2)$. 
\end{claim}

\noindent
\textit{Proof of Claim~\ref{essential}.} 
Suppose for a contradiction that $F$ is compressible. 
Let $D$ be a compressing disk of $F$. 
If $A \cap D = \emptyset$, 
then $D$ is entirely contained in $E(K_i)$ and $F_i$ is compressible, 
contradicting the assumption. 
So in the following we assume $A \cap D \ne \emptyset$. 
Recall that $F \cap E(K_i)$ consists of parallel copies of $F_i$. 
Note that $A  \cap F$ consists of spanning arcs in $A$ 
in minimal number of components (Figure~\ref{arcs}). 
We may assume that $D$ intersects $A$ transversely and the number of components of $A \cap D$ is minimal. 
Then $A \cap D$ consists of circles and arcs whose endpoints belong to $A \cap F$. 
Since $A$ is incompressible, 
we eliminate the circle components, 
and thus $A \cap D$ consists of arcs; 
see Figure~\ref{arcs}. 

\begin{figure}[!ht]
\begin{center}
\includegraphics[width=0.23\linewidth]{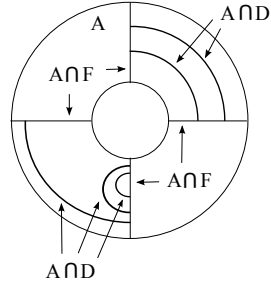}
\caption{$A \cap F$ and $A \cap D$}
\label{arcs}
\end{center}
\end{figure}

Then $A \cap D$ consists of properly embedded arcs in $D$. 
Let $\gamma$ be an outermost arc of $A \cap D$ in $D$; $\gamma$ cuts off an outermost disk $\Delta$. 
There are two possibilities: 
(i) $\partial \gamma$ is contained in a single arc $\tau$ of $A \cap F$ (Figure~\ref{arc_gamma}(i)), 
or (ii) $\gamma$ is an arc connecting two spanning arcs $\tau_1$ and $\tau_2$ of $A \cap F$ (Figure~\ref{arc_gamma}(ii)). 

\begin{figure}[!ht]
\begin{center}
\includegraphics[width=0.7\linewidth]{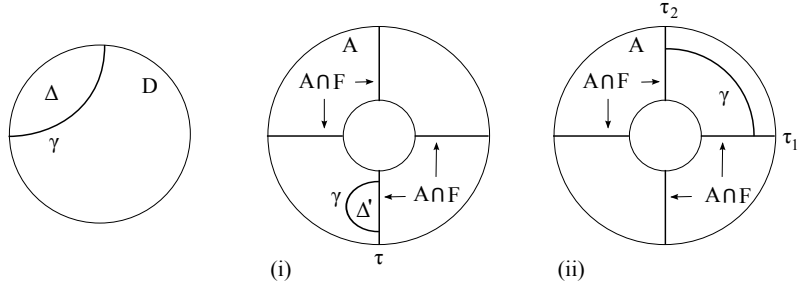}
\caption{An outermost arc $\gamma$ of $A \cap D$ in $D$ and 
its possible situation in $A$}
\label{arc_gamma}
\end{center}
\end{figure}

Suppose that $\partial \gamma$ is contained in a single arc $\tau$ of $A \cap F$ (Figure~\ref{arc_gamma}(i)). 
Then $\gamma$ is parallel to $\tau$; 
$\gamma$ and $\tau$ cobound a disk $\Delta' \subset A$
(If $A \cap F$ consists of a single arc $\tau$, 
then although there would be a possibility that $\gamma$ starts from one side of $\tau$ and ends in the other side 
of $\tau$, this cannot happen for homological reason.)
Let $F_{\Delta}$ be a unique component of $F \cap E(K_1)$ or $F \cap E(K_2)$ intersecting $\partial \Delta$; 
$F_{\Delta}$ is a parallel copy of $F_1$ or $F_2$. 
Then by the incompressibility of $F_i$ in $E(K_i)$ and the irreducibility of $E(K_i)$, 
the disk $\Delta \cup \Delta'$ is parallel to a disk in $\partial E(K_i)$.   
Thus we can isotope $D$ so that $\gamma$ is removed from $A \cap D$. 
This contradicts the minimality of the number of components of $A \cap D$.  

Next assume that $\gamma$ is an arc connecting two spanning arcs $\tau_1$ and $\tau_2$ of $A \cap F$ (Figure~\ref{arc_gamma}(ii)).  
As above we take a unique component $F_{\Delta}$ of $F \cap E(K_1)$ or $F \cap E(K_2)$ intersecting $\partial \Delta$; 
$F_{\Delta}$ is a parallel copy of $F_1$ or $F_2$. 
Since $F_{\Delta}$ is boundary-incompressible, 
$\tau_1$ and $\tau_2$  are contained in a single component of $\partial F_{\Delta}$ and run in opposite directions in $A$, 
a contradiction. 
It follows that $F$ is incompressible in $E(K_1 \sharp K_2)$.  

Since $K_i$ is non-trivial and $F_i$ is not a disk for $i = 1, 2$, 
$F$ is not an annulus. 
Hence \cite[Lemma~1.10]{Hat2}  shows that $F$ is boundary-incompressible as well. 
Thus $F$ is a desired essential surface in $E(K_1 \sharp K_2)$ with boundary slope $p_1/q_1 + p_2/q_2$. 

This completes a proof of Lemma~\ref{boundary slope composite}. 
\hspace*{\fill} $\square$(Lemma~\ref{boundary slope composite})

\medskip

\begin{remark}
\label{nonunique}
In the above construction of the surface $F'$, 
we can slide or twist several times $m_1q_1F_2$ along the annulus $A$ before connecting with $m_2q_2F_1$ 
without changing its boundary slope, 
so $F'$ is not unique. 
\end{remark}

\medskip

Let us turn to a proof of Theorem~\ref{slope conjecture composite}. 
Before proving the theorem, 
we note the following general fact. 

\begin{claim}
\label{mirror}
Let $K$ be a knot in $S^3$. 
If $js^*(K) \subset bs(K)$,  
then $js(K^*) \subset bs(K^*)$. 
\end{claim}

\noindent
\textit{Proof of Claim~\ref{mirror}.}
If $r \in js(K^*)$, 
then $-r \in -js(K^*) = js^*(K) \subset bs(K)$. 
Thus $r \in -bs(K) = bs(K^*)$. 
\hspace*{\fill} $\square$(Claim~\ref{mirror})

\bigskip

\noindent
\textit{Proof of Theorem~\ref{slope conjecture composite}.}
Assume first that $p/q \in js(K_1 \sharp K_2)$. 
Then as shown in Lemma~\ref{jones slope composite}, 
$p/q = p_1/q_1 + p_2/q_2$ for some Jones slopes $p_1/q_1 \in js(K_1)$ and $p_2/q_2 \in js(K_2)$. 
Since $p_1/q_1 \in js(K_1) \subset bs(K_1)$ and $p_2/q_2 \in js(K_2) \subset bs(K_2)$ by the initial assumption, 
Lemma~\ref{boundary slope composite} shows that $p/q = p_1/q_1 + p_2/q_2 \in bs(K_1 \sharp K_2)$. 

Next assume that $p/q \in js^*(K_1 \sharp K_2)$. 
Then $-p/q \in js((K_1 \sharp K_2)^*) = js(K^*_1 \sharp K^*_2)$. 
Since $K_i$ satisfies the slope conjecture, 
$js^*(K_i) \subset bs(K_i)$, thus by Claim~\ref{mirror} $js(K^*_i) \subset bs(K^*_i)$. 
Apply the above argument to $K^*_1$ and $K^*_2$ to conclude that $-p/q \in bs(K^*_1 \sharp K^*_2)= bs((K_1 \sharp K_2)^*)$.  
Hence $p/q \in bs(K_1 \sharp K_2)$. 
This completes a proof of Theorem~\ref{slope conjecture composite}. 
\hspace*{\fill} $\square$(Theorem~\ref{slope conjecture composite})

\section{The slope conjecture and cabling operation}
\label{cable}

Let $V$ be a standardly embedded solid torus in $S^3$ and $k$ a $0$--bridge braid in $V$ which 
wraps $p$ times in meridional direction and $q$ times in longitudinal direction; 
$k$ is a $(p, q)$--torus knot in $S^3$. 
In the following, we assume $q > 1$. 
Given a nontrivial knot $K$, 
take an orientation preserving embedding $f : V \to S^3$ such that the core of $f(V)$ is $K$ and $f$ sends a preferred longitude of $V$ to that of $K$. 
Then the image $f(k)$ is called the \textit{$(p, q)$--cable} of $K$ and denoted by $C_{p, q}(K)$. 
We begin by describing Jones slopes of $C_{p, q}(K)$. 
Let us write $\delta_K(n) = c_2(n)n^2 + c_1(n)n + c_0(n)$, 
where $c_i(n)$ is a periodic function with an integral period. 

\medskip 

In \cite[Proposition~3.2]{KT} Kalfagianni and Tran describe how jones slopes behave under cabling operation. 
See also \cite{KT2}. 
It should be noted here that our normalization of colored Jones functions is slightly different from that in \cite{KT}, 
and $a(n)$, $b(n)$ in \cite{KT} correspond to $c_2(n)$, $c_1(n)+\frac 1 2$, respectively. 

\begin{lemma}[\cite{KT}]
\label{jones slope cable KT} 
Assume that $\delta_K(n)$ has period at most $2$, $c_1(n)+\frac 1 2\le 0$ and 
$4c_2(n) \ne \frac{p}{q}$ for  sufficiently large $n$. 
If $r \in js(C_{p, q}(K))$, 
then $r = pq$ or $aq^2/b$ for some $a/b \in js(K)$. 
\end{lemma}

\medskip

Since the Jones slope $pq$ is the boundary slope of the cabling annulus of $C_{p, q}(K)$, 
$pq \in bs(C_{p, q}(K))$. 
The next result was essentially shown by Klaff and Shalen \cite{KS}, 
but we give a modified proof here.  
See also \cite[Theorem~2.2]{KT}. 

\begin{lemma}
\label{boundary slope cable}
If $a/b \in bs(K)$, 
then $a q^2 / b \in bs(C_{p, q}(K))$. 
\end{lemma}

\noindent
\textit{Proof of Lemma~\ref{boundary slope cable}.}
Let $M_{p, q}  = V - \mathrm{int}N(k)$, the standard $(p, q)$--cable space. 
We denote preferred meridian-longitude pairs of $V$ and $N(k)$ by $(\mu_V, \lambda_V)$ and 
$(\mu, \lambda)$, respectively. 
Then $H_1(M_{p, q}) \cong \mathbb{Z} \oplus \mathbb{Z}$ is generated by $[\lambda_V]$ and $[\mu]$.  
Let $D$ be a $q$--th punctured meridian disk of $V$ and $A$ an obvious annulus connecting $\partial V$ and 
$\partial N(k)$. 
With appropriate orientations we have 
$[D \cap \partial V] = [\mu_V]$, $[D \cap \partial N(k)] = - q [\mu]$,
$[A \cap \partial V] = p[\mu_V] + q[\lambda_V]$, $[A \cap \partial N(k)] = - pq [\mu] - [\lambda]$, 
and thus $[\mu_V] = q[\mu]$, $[\lambda] = q [\lambda_V]$ in $H_1(M_{p, q})$. 

Let $S$ be an oriented surface in $M_{p, q}$ representing the nontrivial homology class 
$(aq-bp)[D] + b[A] \in H_2(M_{p, q}, \partial M_{p, q})$.  
We can construct $S$ by the ``double-curve sum" of $(aq-bp)$ parallel copies of $D$ and $b$ parallel copies of $A$  
(i.e. cut and paste along their intersection arcs to get an embedded surface representing the desired homology class); 
see Figure~\ref{double_curve_sum}. 

\begin{figure}[!ht]
\begin{center}
\includegraphics[width=0.55\linewidth]{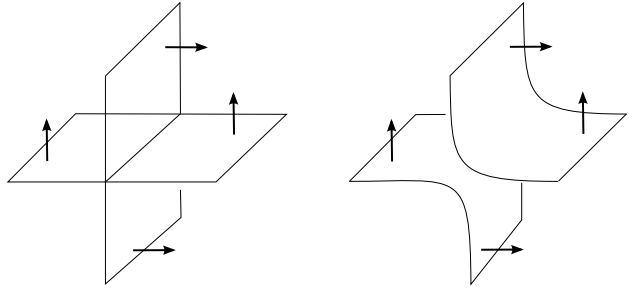}
\caption{Double-curve sum}
\label{double_curve_sum}
\end{center}
\end{figure}

Then it is easy to see that each component of $S \cap \partial V$ has slope $a/b$, 
and that of $S \cap \partial N(k)$ has slope $aq^2/b$. 
If $S$ is compressible, 
then after compression, 
we take a connected component $S_0$ of $S$ which represents nontrivial homology class in 
$H_2(M_{p, q}, \partial M_{p, q})$.  
Since $S_0$ represents a nontrivial homology class, 
it is not a boundary-parallel annulus. 
Thus the incompressible surface $S_0$ is also boundary-incompressible \cite[Lemma~1.10]{Hat2}, 
i.e. $S_0$ is essential in $M_{p, q}$. 
Following \cite[Proposition~1.11]{Hat2}, 
we may assume (up to isotopy) that $S_0$ is horizontal
(i.e. transverse to all Seifert fibers of $M_{p, q}$) or vertical 
(i.e. consists of Seifert fibers of $M_{p, q}$). 
If $S_0$ is vertical, then since $S_0$ is homologically nontrivial, $S_0 = A$. 
If $S_0$ is horizontal, then $S_0 \cap \partial V$ and $S_0 \cap \partial N(k)$ are not empty. 
In particular, each component of $S_0 \cap \partial V$ has slope $a/b$, 
and that of $S_0 \cap \partial N(k)$ has slope $aq^2/b$. 
Let $m_0$ be the number of components of $\partial S_0$ on $\partial V$. 

Now recall that $C_{p, q}(K) = f(k)$, 
where $f : V \to S^3$ is the orientation preserving embedding such that the core of $f(V)$ is $K$. 
Then $E(C_{p, q}(K)) = E(K) \cup f(M_{p, q})$ in which $T = \partial E(K) = \partial f(V)$ is an essential torus. 
Since $a/b \in bs(K)$, 
we have an essential surface $S_1 \subset E(K)$ which has $m_1$ boundary components each of which has slope $a/b$. 
Let us take $m_1$ parallel copies of $f(S_0)$ and $m_0$ parallel copies of $S_1$. 
Connecting them, we obtain a possibly disconnected surface $F'$ in $E(C_{p, q}(K))$. 
Let $F$ be a connected component of $F'$. 
If $F$ is non-orientable, then as in the proof of Lemma~\ref{boundary slope composite} we replace 
$F$ by $\partial N(F)$, where $N(F)$ is a tubular neighborhood of $F$ in $E(C_{p, q}(K))$. 
In the latter case, we continue to use the same symbol $F$ to denote $\partial N(F)$. 
Since $S_0$ and $S_1$ are orientable, 
$F \cap E(K)$ consists of parallel copies of $S_1$, 
and $F \cap f(M_{p, q})$ consists of parallel copies of $f(S_0)$.  
Applying the proof of Claim~\ref{essential}, 
where we use the essentiality of $T$ instead of that of $A$, 
we see that $F$ is an essential surface in $E(C_{p, q}(K))$ with boundary slope $aq^2/b$. 
\hspace*{\fill} $\square$(Lemma~\ref{boundary slope cable})

\medskip

\section{The slope conjecture for graph knots}
\label{graph}

Recall that for any knot $K$, 
$\delta_K(n)$ (resp. $\delta^*_K(n)$) is a quadratic quasi-polynomial $c_2(n)n^2 + c_1(n)n + c_0(n)$ 
(resp. $c^*_2(n)n^2 + c^*_1(n)n + c^*_0(n)$).  

\begin{definition}
\label{condition}
We say that $K$ satisfies \textit{Condition $\delta$} if 
\begin{enumerate}
\item
$\delta_K(n)$ and $\delta^*_K(n)$ have period at most $2$, 
\item
$c_1(n)+\frac 1 2 \le 0$ and $c^*_1(n) -\frac 1 2 \ge 0$, and 
\item $4c_2(n),\ 4c^*_2(n) \in \mathbb{Z}$. 
\end{enumerate}
\end{definition}

\medskip

\begin{proposition}
\label{induction}
Let $\mathcal{K}$ be the maximal set of knots each of which satisfies the slope conjecture and Condition $\delta$. 
Then $\mathcal{K}$ is closed under connected sum and cabling. 
\end{proposition}

\noindent
\textit{Proof of Proposition~\ref{induction}.}
Let us take $K_1, K_2 \in \mathcal{K}$. 
By Theorem~\ref{slope conjecture composite}, $K_1 \sharp K_2$ satisfies the slope conjecture. 
Then it remains to see: 

\begin{claim}
\label{condition sum}
$K_1 \sharp K_2$ satisfies Condition $\delta$. 
\end{claim}

\noindent
\textit{Proof of Claim~\ref{condition sum}.}
Let us write 
$$\delta_{K_1}(n) = \alpha_1(n)n^2 + \beta_1(n)n + \gamma_1(n),$$
$$\delta_{K_2}(n) = \alpha_2(n)n^2 + \beta_2(n)n + \gamma_2(n).$$
Then we have: 
$$\delta_{K_1 \sharp K_2}(n) = (\alpha_1(n)+\alpha_2(n))n^2+(\beta_1(n)+\beta_2(n))n+(\gamma_1(n)+\gamma_2(n)).$$ 
Since the common period of $\alpha_i(n), \beta_i(n)$ and $\gamma_i(n)$ is at most $2$, 
$\alpha_1(n)+\alpha_2(n),\ \beta_1(n)+\beta_2(n)$ and $\gamma_1(n)+\gamma_2(n)$ have period at most $2$, 
and hence $\delta_{K_1 \sharp K_2}(n)$ has also period $\le 2$. 
Since $\beta_1(n) + \frac 1 2\le 0$ and $\beta_2(n)+\frac 1 2 \le 0$, 
$(\beta_1(n) + \beta_2(n)) + \frac 1 2 \le 0$, which shows $(2)$. 
It is obvious that $4(\alpha_1(n) + \alpha_2(n))$ is an integer. 
It is easy to check the remaining conditions in a similar fashion. 
\hspace*{\fill} $\square$(Claim~\ref{condition sum})  

\medskip

\begin{claim}
\label{condition cable}
Let us take $K \in \mathcal{K}$ and its cable $C(K)$. 
Then $C(K) \in \mathcal{K}$. 
\end{claim}

\noindent
\textit{Proof of Claim~\ref{condition cable}.}
This was shown in \cite{KT}. 
Combining Lemmas~\ref{jones slope cable KT}, \ref{boundary slope cable} and Claim~\ref{mirror}  
we see that $C(K)$ satisfies the slope conjecture \cite[Theorem~3.4]{KT}. 
Proposition~3.2 in \cite{KT}, together with mirroring technique, 
shows that $C(K)$ satisfies Condition $\delta$ as well. 
\hspace*{\fill} $\square$(Claim~\ref{condition cable})  

\medskip

This establishes Proposition~\ref{induction}.
\hspace*{\fill} $\square$(Proposition~\ref{induction})

\medskip

\noindent
\textit{Proof of Theorem~\ref{slope conjecture graph}.}
Let $K$ be a graph knot.
If $K$ is the trivial knot,  $K$ obviously satisfies the slope conjecture 
($js(K) \cup js^*(K) = bs(K) = \{ 0 \}$). 
If $K$ is nontrivial, 
then $K$ is obtained from torus knots by a finite sequence of operations of cabling 
and connected sum; see \cite[Corollary~4.2]{Gor}. 
Garoufalidis \cite[4.8]{Garoufalidis} proves the slope conjecture for torus knots. 
Actually he computes their colored Jones functions of $T_{p, q}$ $(p, q > 0)$ explicitly: 
\begin{eqnarray*}
 && \delta_{T_{p,q}}(n) = \frac{pq}{4}n^2 - \frac{1}{2}n - \frac{pq-2}{4} - (1+(-1)^n)\frac{(p-2)(q-2)}{8},\\
&& \delta^*_{T_{p, q}}(n) = \frac{(p-1)(q-1)}{2}n - \frac{(p-1)(q-1)}{2}. 
\end{eqnarray*}
Then it is easy to see that $T_{p, q}$ satisfies Condition $\delta$. 
Since $\delta_{T_{-p, q}}(n) = - \delta^*_{T_{p, q}}(n)$ and $\delta^*_{T_{-p, q}}(n) = - \delta_{T_{p, q}}(n)$, 
any nontrivial torus knot satisfies Condition $\delta$. 
It follows from Proposition~\ref{induction} that the set of nontrivial graph knots is contained in $\mathcal{K}$. 
Thus any graph knot satisfies the slope conjecture. 
\hspace*{\fill} $\square$(Theorem~\ref{slope conjecture graph})

\end{document}